\documentclass[letterpaper, 10pt, conference]{ieeeconf}
\IEEEoverridecommandlockouts 
\usepackage[noformatting,nogeometry]{stdcommands}
\usepackage{graphicx}
\usepackage{dsfont} 
\usepackage{amsmath}
\usepackage{amssymb}
\usepackage{url}
\usepackage[utf8]{inputenc} 
\usepackage[format=hang,singlelinecheck=false,caption=false,font=footnotesize]{subfig}
\usepackage{color}
\usepackage{algorithm}
\usepackage{algpseudocode}

\graphicspath{{./images/}}

\renewcommand{\one}{\mathds{1}}

\newtheorem{thm}{Theorem}[section]

\newtheorem{definition}[thm]{Definition}
\newtheorem{assumption}[thm]{Assumption}
\newtheorem{prop}[thm]{Proposition}

\DeclareMathOperator{\supp}{supp}
\DeclareMathOperator{\cohull}{co}

\title{\LARGE \bf
  Numerical Synthesis of Pontryagin Optimal Control Minimizers Using Sampling-Based Methods
}

\author{
  Runxin He \and
  Humberto Gonzalez%
  \thanks{
    The work was supported in part by the U.S.\ National Science Foundation award \#1646579.
    Runxin He (\texttt{r.he@wustl.edu}) is with the Department of Electrical \& Systems Engineering, Washington University in St.\,Louis, St.\,Louis, MO 63130, USA.
    Humberto Gonzalez is a Senior Data Scientist at Instart Logic, Palo Alto, CA 94306, USA.
  }%
}

\begin{document}

\maketitle
\thispagestyle{empty}
\pagestyle{empty}
\begin{abstract}
  We present a theoretical formulation, and a corresponding numerical algorithm, that can find Pontryagin-optimal inputs for general dynamical systems by using a direct method.
  Optimal control remains as a versatile and relevant framework in systems theory applications, many decades after being formally defined.
  Pontryagin-optimal inputs can be found for some classes of problems using indirect methods, but these are often slow or lack robustness.
  On the other hand, convergent direct optimal control methods are fast, but their solutions usually converge to first-order optimality conditions, which are weaker.

  Our result, founded on the theory of relaxed inputs as defined by J.\,Warga, establishes an equivalence between Pontryagin-optimal inputs and optimal relaxed inputs.
  We also formulate a sampling-based numerical method to approximate the Pontryagin-optimal relaxed inputs using an iterative direct method.
  Finally, using a provably-convergent numerical method, we synthesize approximations of the Pontryagin-optimal inputs from the sampled relaxed inputs.
\end{abstract}

\section{INTRODUCTION}
\label{sec:introduction}

Optimal control is a theoretical and practical framework that has been widely used to analyze the behavior of controlled dynamical systems~\cite{zhou1996robust, bertsekas1995dynamic}, and to synthesize actuation actions for dynamical systems in the face of safety, robustness, or uncertainty considerations~\cite{bryson1975applied, bellman2015applied}.
While some optimal control problems can be solved using purely analytical or algebraic tools~\cite{boyd1991linear}, modern computers and dynamical systems embedded in changing environments have led to a surge in numerical methods for optimal control~\cite{rao2009survey}.

Numerical methods for optimal control are typically divided into indirect methods, where the optimal solution is found as the solution of a set of equations typically derived from necessary optimality conditions, and direct methods, where the solution is found by iteratively minimizing the cost function at hand.
As shown by Polak~\cite[Sec.~4.2.6]{polak2012optimization} and Schwartz~\cite{schwartz1996theory}, direct numerical methods based on explicit time discretization converge to solutions satisfying derivative necessary optimality condition, which are strictly weaker than Minimum Principle condition.
Although some direct methods result in points satisfying the Minimum Principle, their implementation is usually impractical~\cite{gonzalez2010numerical}.
Some indirect methods also converge to inputs satisfying the Minimum Principle, usually relying on multiple-point boundary-value problems~\cite{betts2010practical}, which only converge for suitably chosen initial guesses~\cite{betts2003initialization}.

This paper presents a theoretical and numerical framework, which results in an algorithm that uses the direct method to converge to Pontryagin Minimum Principle inputs.
Our result is founded on the theory of relaxed control~\cite{warga2014optimal, warga1962relaxed}, which we use to derive a convergent sampling-based numerical method.
Once a relaxed control has been numerically computed, we use a projection operation originally devised for switched dynamical systems~\cite{vasudevan2013a, vasudevan2013b} to synthesize arbitrarily accurate approximations of the trajectories generated by the relaxed inputs.
Our result bridges a significant gap between the formulation of conceptual algorithms and implementable numerical algorithms which converge to Minimum Principle through direct method.
Moreover, our method achieves its goals in a numerically efficient and scalable way.

Our paper is organized as follows.
Sec.~\ref{sec:problem} introduces our notation and the optimal control problem we aim to solve.
A conceptual algorithm is presented in Sec.~\ref{sec:method}, and Sec.~\ref{sec:implementation} describes an implementable sampling-based numerical algorithm.
Finally, simulation results are shown in Sec.~\ref{sec:examples}.
Due to space limitations we omit all the formal proofs in this paper.

\section{PROBLEM DESCRIPTION}
\label{sec:problem}


We begin by introducing notation and preliminary results necessary to formulate our theoretical and numerical results in Sec.~\ref{sec:prelim}.
Then, in Sec.~\ref{sec:conceptual}, we formulate the conceptual optimal control problem we will address throughout this paper.

\subsection{Preliminaries}
\label{sec:prelim}

Given $n \in \N$ and $p \geq 1$, we denote the standard finite-dimensional $p$-norm by $\norm{\cdot}_p$, and the induced matrix $p$-norm by $\norm{\cdot}_{i,p}$.
We denote by ${\cal M}(\R^n)$ the set of Radon measures defined over the Borel sets of $\R^n$.
Given $\mu \in {\cal M}(\R^m)$ for $m \in \N$, we say that a function $f\colon \R^m \to \R^n$ is $L^2_\mu$-integrable, denoted $f \in L^2_\mu\p{\R^m, \R^n}$, if there exists $p \geq 1$ such that  $\norm{f}_{\mu} = \p{\int_{\R^m} \norm{f(x)}_p^2\, \diff{\mu}(x)}^{\!1/2}$ is finite.
We abuse notation and denote by $L^2\p{\R^m,\R^n}$ the space of Lebesgue square-integrable functions.
Furthermore, we say that $\mu \in {\cal M}(\R^m)$ is a probability measure if $\mu(\R^m) = 1$.
We denote the set of all probability measures by ${\cal M}_p(\R^m)$.
A stochastic process is a function $\mu\colon [0,T] \to {\cal M}_p(\R^m)$, and throughout the paper we simply write $\mu_t$ instead of $\mu(t)$.


Let $\mu_1,\mu_2 \in {\cal M}(\R^m)$ be two Radon measures.
Then the difference $\nu = \mu_1 - \mu_2$ is a signed measure, and we define $L^2_{\nu}(\R^m, \R^n) = L^2_{\mu_1}\p{\R^m,\R^n} \cap L^2_{\mu_2}\p{\R^m,\R^n}$.
Given $f \in L^2_{\nu}\p{\R^m,\R^n}$, its integral with respect to $\nu$ is defined by
$\int_{\R^m} f(x)\, \diff{\nu}(x) = \int_{\R^m} f(x)\, \diff{\mu_1}(x) - \int_{\R^m} f(x)\, \diff{\mu_2}(x)$.


\subsection{Conceptual Optimal Control Problem}
\label{sec:conceptual}

We are interested in solving an input-constrained optimal control problem, formulated as follows:
\begin{equation}
  \label{eq:original}
  \begin{aligned}
    \min_{\substack{%
        x \in L^2([0,T],\R^n)\\
        u \in L^2([0,T],\R^m)
      }}\,
    & \Psi\pb{x(T)},\\
    \text{subject to:}\quad
    & x(0) = \xi,\\
    & \dot{x}(t) = f\pb{t, x(t), u(t)},\\
    & u(t) \in U,\ \text{for a.\,e.}\ t \in [0,T],
  \end{aligned}
\end{equation}
where $U$ is a connected and compact subset of $\R^m$, $\xi \in \R^n$, and the functions $f$ and $\Psi$ are well-defined, each with an appropriate domain and range.
We will say that $f$ is the vector field and $\Psi$ the final cost of the problem in eq.~\eqref{eq:original}.
Note that the optimal control problem in eq.~\eqref{eq:original} is quite general, since other standard formulations, such as those including running cost functions or minimum-time cost functions, can be easily converted to it (see~\cite[Sec.~4.1.2]{polak2012optimization}).

First we give the following assumption to guarantee the uniqueness of the trajectories, as well as the convergence of our numerical method, in this paper.
\begin{assumption}
  \label{assump:lipschitz}
  The functions $f$ and $\Psi$ are Lipschitz continuously differentiable.
  Thus, there exists $L > 0$ such that, for each $t_1, t_2 \in [0,T]$, $x_1, x_2 \in \R^n$, and $u_1, u_2 \in U$:
    \begin{align}
    &\abs{\Psi(x_1) - \Psi(x_2)} \leq L\, \norm{x_1 - x_2}_2,\\
    &\normb{\spderiv{\Psi}{x}(x_1) - \spderiv{\Psi}{x}(x_2)}_2 \leq L\, \norm{x_1 - x_2}_2,\\
    &\norm{f(t_1, x_1, u_1) - f(t_2, x_2, u_2)}_2 \notag\\
    &\hspace{.5in} \leq L\, \pb{\abs{t_1 - t_2} + \norm{x_1 - x_2}_2 + \norm{u_1 - u_2}_2},\\
    &\normb{\spderiv{f}{x}(t_1, x_1, u_1) - \spderiv{f}{x}(t_2, x_2, u_2)}_{i,2} \notag\\
    &\hspace{.5in} \leq L\, \pb{\abs{t_1 - t_2} + \norm{x_1 - x_2}_2 + \norm{u_1 - u_2}_2},\\
    &\normb{\spderiv{f}{u}(t_1, x_1, u_1) - \spderiv{f}{u}(t_2, x_2, u_2)}_{i,2} \notag\\
    &\hspace{.5in} \leq L\, \pb{\abs{t_1 - t_2} + \norm{x_1 - x_2}_2 + \norm{u_1 - u_2}_2}.
  \end{align}
\end{assumption}

Now, we relax the problem in eq.~\eqref{eq:original} using the concept of relaxed inputs.
Consider the following relaxed optimal control problem:
\begin{equation}
  \label{eq:original2}
  \begin{aligned}
    \min_{\substack{%
        x \in L^2([0,T],\R^n)\\
        \mu\colon [0,T] \to {\cal M}_p(\R^m)
      }}\,
    & \Psi\pb{x(T)},\\
    \text{subject to:}\quad
    & x(0) = \xi,\\
    & \dot{x}(t) = \int_{\mathrlap{\R^m}}\ f\pb{t, x(t), u}\, \diff{\mu_t}(u),\\
    & \supp(\mu_t) \subset U,\ \text{for a.\,e.}\ t \in [0,T],
  \end{aligned}
\end{equation}
where $\supp(\mu_t)$ is the support of $\mu_t$, i.e., the smallest set $S$ such that $\mu_t(S) = 1$.
In other words, instead of optimizing over the space of $L^2$ functions, we optimize over the space of stochastic processes defined on $U$.

As shown in~\cite[Theorem II.6.5]{warga2014optimal}, if the vector field $f(t, x, u)$ satisfies Asm.~\ref{assump:lipschitz}, the relaxed problem in eq.~\eqref{eq:original2} always has a
solution.
Given a fixed initial condition $x(0) = \xi$, we denote the unique trajectory resulting from the stochastic process $\mu$ by $x^{(\mu)}$.
For simplicity, we also denote the unique trajectory resulting from an input $u$ by $x^{(u)}$.

Note that the problem in eq.~\eqref{eq:original} is a particular case of the problem in eq.~\eqref{eq:original2}.
Indeed, given an arbitrary input $\hat{u} \in L^2([0,T],\R^m)$, the stochastic process defined by $\mu_t(S) = 1$ whenever $\hat{u}(t) \in S$, and $\mu_t(S) = 0$ otherwise, produces exactly the same trajectory
as $\hat{u}$.
This implies that the feasible set of the relaxed problem is strictly larger than that of the original problem.
However, both original and relaxed problems result in the same optimal values.
Their equivalence follows since every point in the feasible set of the relaxed problem can be arbitrarily approximated using points in the feasible set of the original problem.
The following theorem is an extension of the Chattering Lemma~\cite[Thm.~4.1]{berkovitz2013optimal}.
\begin{thm}
  \label{thm:approx}
  Let $f$ be a vector field satisfying assumption~\ref{assump:lipschitz}, and let $\mu\colon [0,T] \to {\cal M}_p(\R^m)$ be a stochastic process.
  Then, for each $\epsilon > 0$ there exists a control signal $\tilde{u}(t)$ such that for each $t \in [0,T]$,
  $\normb{x^{(\mu)}(t) - x^{(\tilde{u})}(t)}_2 < \epsilon$.
\end{thm}

\section{OPTIMALITY CONDITIONS FOR OPTIMAL CONTROL}
\label{sec:method}


As shown by the Minimum Principle~\cite{pontryagin1987mathematical}, it is possible to find a necessary condition for optimal points that cannot be formulated using directional derivatives, in contrast of variational-based necessary conditions for optimal points that can be viewed as extensions of finite-dimensional first-order (or KKT) conditions~\cite{nocedal2006}.
Furthermore, necessary conditions for the problem in eq.~\eqref{eq:original} based on directional derivatives are strictly weaker than those based on the Minimum Principle.
Indeed, consider the following optimal control problem:
\begin{equation}
  \label{eq:example}
  \min\set{x(1) \mid
    x(0) = 0,\
    \dot{x}(t) = \sfrac{1}{2}\, \abs{u(t)} - \cos u(t)}.
\end{equation}
Note that the costate associated to $x$ is $p(t) = 1$ for each $t$, hence the Hamiltonian of this system is identical to its vector field, as defined in eqs.~\eqref{eq:hamiltonian0} and~\eqref{eq:adjoint0} respectively.
In this example, the Minimum Principle results in a single minimizer $u_g(t) = 0$, while first-order methods will converge to other local minimizers of the form $u_l(t) \in \set{2\, n\, \pi - \sfrac{\pi}{6}, \sfrac{\pi}{6} - 2\, n\, \pi \mid n = 1,2,\dotsc}$ depending on the initialization of the optimization algorithm.

In this section we present the theoretical foundation for our numerical method, including a conceptual infinite-dimensional optimization algorithm.


\subsection{Optimality Functions}
\label{sec:optfunc}

\begin{definition}[Sec.~1.2 in \cite{polak2012optimization}]
  Consider an optimization problem defined on ${\cal X}$.
  Then $\theta\colon {\cal X} \to (-\infty,0]$ is an optimality function iff $\theta(\hat{x}) = 0$ for each minimizer $\hat{x} \in {\cal X}$.
\end{definition}

Optimality functions are useful in practice since $\theta(x) < 0$ implies $x$ is \emph{not} a minimizer.
Hence, they can be used as numerical tests to check whether a minimizer has been reached.

Given an input $u_0$ and its trajectory $x^{(u_0)}$, consider:
\begin{equation}
  \begin{aligned}
    \theta_{o, l}\pb{x^{(u_0)}, u_0}
    = &\min_{\substack{%
        \delta x \in L^2([0,T],\R^n)\\
        \delta u \in L^2([0,T],\R^m)
      }}\,
    \spderiv{\Psi}{x}\pb{x^{(u_0)}(T)}^\tp\, \delta x(T) \\
    \text{subject to:} \quad
    &\delta x (0) = 0, \\
    &\delta \dot{x}(t) = \spderiv{f}{x}\pb{t, x^{(u_0)}(t), u_0(t)}\, \delta x(t)\\
    &\phantom{\delta \dot{x}(t) =} + \spderiv{f}{u}\pb{t, x^{(u_0)}(t), u_0(t)}\, \delta u(t), \\
    &u_0(t) + \delta u(t) \in U,\ \text{for a.e.}\ t \in [0,T],
  \end{aligned}
  \label{eq:thetaod}
\end{equation}
and:
\begin{equation}
  \label{eq:thetaop}
  \begin{aligned}
    \theta_{o, h}\pb{x^{(u_0)}, u_0}
    = &\hspace{-6pt}\min_{u \in L^2([0,T],\R^m)}
    \int_{0}^{\mathrlap{T}}\!\! H_0\pb{t, x^{(u_0)}(t), u(t), p_0(t)}\\
    &\hspace{3em} - H_0\pb{t, x^{(u_0)}(t), u_0(t), p_0(t)}\, \diff{t},\\
    \text{subject to:} \quad
    & u(t) \in U,\ \text{for a.\,e.}\ t \in [0,T],
  \end{aligned}
\end{equation}
where $H_0$ is the Hamiltonian of problem~\eqref{eq:original} at time $t$:
\begin{equation}
  \label{eq:hamiltonian0}
  H_0\pb{t, x^{(u)}(t), u(t), p_0(t)}
  = p_0(t)^\tp\, f\pb{t, x^{(u)}(t), u(t)},
\end{equation}
and $p_0(t)$ is the costate of problem~\eqref{eq:original}, defined by:
\begin{equation}
  \label{eq:adjoint0}
  \begin{aligned}
    p_0(T) &= \spderiv{\Psi}{x}\pb{x^{\p{u_0}}(T)},\\
    \dot{p}_0(t) &= - \spderiv{f}{x} \pb{t, x^{\p{u_0}} (t), u_0 (t)}^\tp\, p_0(t).
  \end{aligned}
\end{equation}

We omit the proofs of the following propositions, but they follow closely the arguments in Thms.~5.6.8 and~5.6.9 in~\cite{polak2012optimization}, and Prop.~4.5 in~\cite{Polak1984}.
\begin{prop}
  \label{prop:theta_original}
  The functions $\theta_{o,l}$ and $\theta_{o,h}$, defined in eqs.~\eqref{eq:thetaod} and~\eqref{eq:thetaop}, are optimality functions of the problem~\eqref{eq:original}.
\end{prop}
\begin{prop}
  \label{prop:thetao_equiv}
  Let $u_0$ be an input and $x^{(u_0)}$ its trajectory.
  If $\theta_{o, h}\pb{x^{(u_0)}, u_0}$ in eq.~\eqref{eq:thetaop} equals zero, then $\theta_{o, l} \p{x^{\p{u_0}}, u_{0}}$ in eq.~\eqref{eq:thetaod} also equals zero.
\end{prop}

Note that the opposite statement to Prop.~\ref{prop:thetao_equiv} is not true in general~\cite[Sec.~4.2.6]{polak2012optimization}.
This is a significant practical problem, since $\theta_{o,l}$, which captures first-order minimizers, is a convex problem, while $\theta_{o,h}$, which captures the stronger Pontryagin minimizers, is in general a non-convex problem.
Thus, direct numerical methods to compute the optimal control of the original formulation in~\eqref{eq:original} can only provably capture Pontryagin minimizers when the non-convex problem in eq.~\eqref{eq:thetaop} can be simplified or solved analytically.

Now, let us define an optimality function for the relaxed problem in equation~\eqref{eq:original2}.
Given a stochastic process $\mu_{0}$ and its corresponding trajectory $x^{\p{\mu_{0}}}$, the Hamiltonian of problem~\eqref{eq:original2} at time $t$ is:
\begin{equation}
  \label{eq:hamiltonian}
  H\pb{t, x^{(\mu)}(t), \mu_t, p(t)}
  = p(t)^\tp \int_{\mathrlap{\R^m}}\, f\pb{t, x^{(\mu)}(t), u}\, \diff{\mu_t}(u),
\end{equation}
where $p(t)$ is the costate defined by:
\begin{equation}
  \label{eq:adjoint}
  \begin{aligned}
    p(T) &= \spderiv{\Psi}{x}\pb{x^{(\mu)}}(T),\\
    \dot{p}(t) &= - \int_{\R^m}\! \spderiv{f}{x}\pb{t, x^{(\mu)}(t), u}^\tp\, \diff{\mu_t}(u)\, p(t).
  \end{aligned}
\end{equation}

Given a stochastic process $\mu_0$ and its trajectory $x^{(\mu_0)}$, consider:
\begin{equation}
  \label{eq:lin_optimal_mu}
  \begin{aligned}
    \theta_l\pb{x^{(\mu_0)},\mu_0} =
    &\min_{\substack{%
        \delta x \in L^2([0,T],\R^n)\\
        \delta \mu\colon [0,T] \to {\cal M}(\R^m)
      }}\,
    \spderiv{\Psi}{x}\pb{x^{(\mu_0)}(T)}^\tp\, \delta x(T)
    \\
    \text{subject to:}\quad
    &\delta\dot{x}(t) =
    \int_{\mathrlap{\R^m}}\ \, \spderiv{f}{x}\pb{t, x^{(\mu_0)}(t), u}\, \diff{\mu_{0, t}}(u)\, \delta x(t)\\
    &\phantom{\dot{\delta x} \p{t} =}\
    + \int_{\mathrlap{\R^m}}\, f\pb{t, x^{\p{\mu_0}} \p{t}, u}\, \diff{\delta \mu_t}(u), \\
    &\delta x \p{0} = 0, \\
    &\supp(\delta \mu_{t} + \mu_{0, t}) \subset U,\\
    &\delta\mu_t + \mu_{0,t} \geq 0,\\
    &\int_{\mathrlap{\R^m}}\ \diff{\delta\mu_t}(u) = 0\ \text{for a.\,e.}\ t \in [0,T],
  \end{aligned}
\end{equation}
and:
\begin{equation}
  \label{eq:ham_optimal_mu}
  \begin{aligned}
    \theta_{h}\pb{x^{\p{\mu_0}}, \mu_{0}}\!
    =&\hspace{-5pt}
    \min_{\delta \mu\colon\! [0,T] \to {\cal M}(\R^m)}\!
    \int_0^{\mathrlap{T}} \!\!H\pb{t, x^{(\mu_0)}(t), \delta\mu_t, p(t)} \diff{t},\\
    \text{subject to:} \quad
    & \supp(\delta \mu_{t} + \mu_{0, t}) \subset U,\\
    & \delta\mu_t + \mu_{0,t} \geq 0,\\
    & \int_{\mathrlap{\R^m}}\ \diff{\delta\mu_t}(u) = 0,\ \text{for a.e.}\ t \in [0,T],
  \end{aligned}
\end{equation}
where $\delta \mu_t$ is a signed measure in $\R^m$.

\begin{prop}
  \label{prop:optimality_hamiltonian}
  The functions $\theta_l$ and $\theta_h$, defined in eqs.~\eqref{eq:lin_optimal_mu} and~\eqref{eq:ham_optimal_mu}, are optimality functions of the problem~\eqref{eq:original2}.
\end{prop}

Similar to the optimality functions for problem~\eqref{eq:original}, $\theta_l$ in eq.~\eqref{eq:lin_optimal_mu} captures first-order minimizers, while $\theta_h$ in eq.~\eqref{eq:ham_optimal_mu} captures Pontryagin minimizers.
In the next subsection we argue that $\theta_h$ is in fact equivalent to $\theta_l$, thus either can be used in practical and implementable numerical algorithms.

\subsection{Gradient Descent Method for Relaxed Problems}
\label{sec:graddesc}

Now we can show the connection between the directional and Pontryagin optimality functions.
\begin{thm}
  \label{thm:eqv_optimal}
  The optimality functions $\theta_l$ in eq.~\eqref{eq:lin_optimal_mu} and $\theta_h$ in eq.~\eqref{eq:ham_optimal_mu} are equivalent.
  That is, given a pair $\pb{x^{(\mu)},\mu}$, both optimality functions produce the same value and minimizers.
\end{thm}

We omit a detailed proof.
Nonetheless, the argument follows using the costate in eq.~\eqref{eq:adjoint} to derive the Fr\'echet derivative of $\Psi$ as in~\cite[Thm.~5.6.9]{polak2012optimization}, and then rewrite $\theta_l$ using the costate.

A significant feature of optimality functions based on first-order derivatives is that their minimizing arguments are also descent directions for the objective function.
The following proposition shows that this property is preserved by $\theta_l$.
\begin{prop}
  Let $\mu_0$ be a stochastic process and $x^{(\mu_0)}$ its corresponding trajectory.
  Suppose that $\pb{\delta x^*, \delta \mu_0^*}$ is the minimizing argument of $\theta_l$ in~\eqref{eq:lin_optimal_mu}.
  Then there exists a step size $\lambda \in (0, 1)$ such that:
  \begin{equation}
    \Psi\pb{x^{(\mu_0 + \lambda\, \delta \mu_0^*)}(T)} \leq \Psi\pb{x^{(\mu_0)}(T)}.
  \end{equation}
\end{prop}

In practice, given $\alpha, \beta \in \p{0, 1}$, the step size $\lambda$ can be obtained using the following variation of the Armijo algorithm~\cite{Armijo1966minimization}:
\begin{multline}
  \label{eq:arimijo_mu}
  \lambda^* = \min\setb{\beta^k \mid
    \Psi\pb{x^{(\mu_{0} + \beta^k\, \delta \mu_0^*)}(T)} - \Psi\pb{x^{(\mu_0)}(T)}\\
    \leq \alpha\, \beta^k\, \theta_h\pb{x^{(\mu_0)}, \mu_0},\ k \in \N},
\end{multline}

\begin{figure}[tp]
  \begin{algorithmic}[1]
    \Require $\mu\colon [0,T] \to {\cal M}_p(\R^m)$, $\alpha, \beta \in (0, 1)$.
    \Loop
    \State Compute $x^{(\mu)}$.
    \State Compute $\theta_h\pb{x^{(\mu)}, \mu}$ and $\delta \mu$ as in~\eqref{eq:ham_optimal_mu}.
    \If{$\theta_h\pb{x^{(\mu)}, \mu} = 0$}
    \State \Return $\mu$.
    \EndIf
    \State Compute $\lambda^*$ as in~\eqref{eq:arimijo_mu}.
    \State $\mu \gets \mu + \lambda^*\, \delta \mu$.
    \EndLoop
  \end{algorithmic}
  \caption{Conceptual optimization algorithm to solve problem~\eqref{eq:original2}.}
  \label{alg:infinite}
\end{figure}

Fig.~\ref{alg:infinite} shows a conceptual algorithm to solve the relaxed optimal control problem in~\eqref{eq:original2}.
The theorem below shows its convergence, whose proof follows thanks to Thm.~\ref{thm:eqv_optimal} and the argument in the proof of~\cite[Thm.~5.12]{vasudevan2013a}.

\begin{thm}
  \label{thm:convegence}
  Let $\set{\mu_i}_{i \in \N}$ be a sequence of stochastic processes generated by the algorithm in Fig.~\ref{alg:infinite}, and let $\setn{x^{(\mu_i)}}_{i \in \N}$ be its corresponding sequence of trajectories.
  Then $\lim_{i \to \infty} \theta_h\pb{x^{(\mu_i)}, \mu_i} = 0$.
\end{thm}

\section{SYNTHESIS OF RELAXED OPTIMAL INPUTS}
\label{sec:implementation}

Now that we have established a theoretical foundation for the equivalence between first-order derivative and Minimum Principle minimizers for relaxed optimal control problems, we focus our attention on the development of numerical algorithms to synthesize approximated optimal control inputs.
The iteration algorithm in this section is based on solving a sequence of convex optimization problems, even when the dynamical system is nonlinear.

As we show in Fig.~\ref{alg:infinite}, it is possible to formulate an iterative gradient descent method using relaxed inputs that converges to Pontryagin minimizers.
This conceptual algorithm requires solving the optimization problem $\theta_h$, which generates a variation of a stochastic process to locally reduce the cost function of the problem in eq.~\eqref{eq:original2}.

The following propositions enable us to find efficient numerical implementations of our algorithm in Fig.~\ref{alg:infinite}.
\begin{prop}
  \label{prop:theta_convex}
  $\theta_h$ is a convex optimization problem.
\end{prop}
\begin{prop}
  \label{prop:cohull}
  For each $x \in \R^n$:
  \begin{multline}
     \set{\int_{\mathrlap{\R^m}}\, f(x,u)\, \diff{\mu}(u) \mid \mu \in {\cal M}(\R^m),\ \supp(\mu) \subset U}\\
    = \cohull\set{f(x,u) \mid u \in U},
  \end{multline}
  where $\cohull(S)$ is the convex hull of $S \subset \R^n$.
\end{prop}

Using Prop.~\ref{prop:cohull} and eq.~\eqref{eq:ham_optimal_mu} we have almost all the necessary results to develop an implementable version of the algorithm in Fig.~\ref{alg:infinite}.
We now focus our attention on the two remaining problems: how to approximate the convex hull in Prop.~\ref{prop:cohull}, and how to synthesize control inputs given a stochastic process $\mu_t$.



\subsection{Convex Hull Approximation}
\label{subsec:sampling}

\begin{figure}[tp]
  \centering
  \scalebox{0.8}{\input{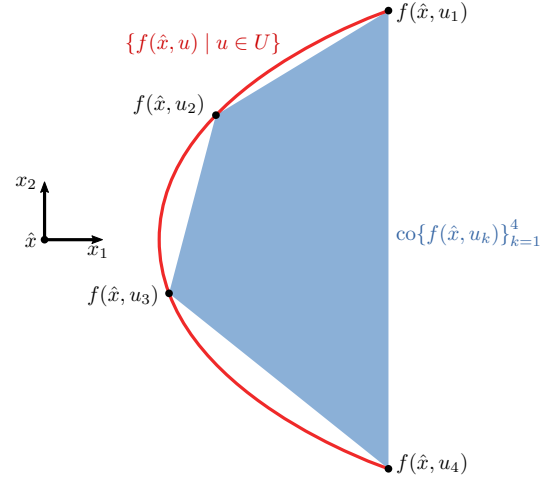}}%
  \caption{%
    Vector field set at $\hat{x}$ (red), and its sampled convex hull (blue), for $f(x,u) = x + \p{u^2 + 1, u}^\tp$, $U = \sp{-1,1}$.
  }
  \label{fig:convex_hull}
\end{figure}

An effective method to find an inner approximation of a convex set is using the convex hull of samples obtained from the set, as shown in Fig.~\ref{fig:convex_hull}.
In this paper we use a simple Monte Carlo sampling method.
Note that other sampling methods, such as Importance Sampling or MCMC~\cite{andrieu2003introduction}, can also be used and might result in faster computations for large values of the input dimension $m$.

Given an $\mu$ and its trajectory $x^{(\mu)}$, we uniformly obtain $N$ samples from the input set $\set{u_i}_{i=1}^N \subset U$ to create a set of vector field evaluations $\set{f\pb{t, x^{(\mu)}(t), u_i}}_{i=1}^N$.
We build an approximated stochastic process $\mu_N$ from these samples using a set of weight functions $\set{w_i(t)}_{i=1}^N$ such that $\sum_{i=1}^N w_i(t) = 1$ and $w_i(t) \geq 0$ for each $i \in \set{1, \dotsc, N}$.

\begin{prop}
  \label{prop:sampling_converge}
  Let $\mu \in {\cal M}_p(\R^m)$ be a stochastic process with $\supp(\mu) \subset U$, and let $\set{u_i}_{i=1}^N$ be a set of samples uniformly drawn from $U$.
  Then there exists $\set{w_{N,i}(t)}_{i=1}^N$ such that the empirical stochastic process $\mu_{N,t} = \sum_{i=1}^N w_{N,i}(t)\, \one_{u_i} \in {\cal M}_p(\R^m)$ satisfies:
  \begin{equation}
    \sum_{i = 1}^N w_{N,i}(t)\, f\pb{t, x^{(\mu_N)}(t), u_i}
    \to
    \int_{\mathrlap{\R^m}}\, f\pb{t, x^{(\mu)}(t), u}\, \diff{\mu_t}(u),
  \end{equation}
  as $N \to \infty$ for a.e.\ $t \in [0,T]$, where $\one_u$ is the Dirac delta function at $u \in U$.
\end{prop}

Consider the following discretized version of the optimality function $\theta_h$ in eq.~\eqref{eq:ham_optimal_mu}:
\begin{equation}
  \label{eq:ham_optimal_muN}
  \begin{aligned}
    \theta_{N,h}\pb{&x^{(\mu_N)}, \mu_N}\!\\
    &= \min_{\delta\mu_N\colon\! [0,T] \to {\cal M}(\R^m)}
    \int_0^{\mathrlap{T}} H\pb{t, x^{(\mu_N)}(t), \delta\mu_{N,t}, p(t)} \diff{t},\\
    &\text{subject to:}\quad
    \supp(\delta \mu_{N,t}) \subset \supp(\mu_{N, t}),\\
    &\phantom{\text{subject to:}\quad}
    \delta\mu_{N,t} + \mu_{N,t} \geq 0,\\
    &\phantom{\text{subject to:}\quad}
    \int_{\mathrlap{\R^m}}\ \diff{\delta\mu_{N,t}}(u) = 0,\ \text{for a.e.}\ t \in [0,T],
  \end{aligned}
\end{equation}
We can now establish the epi-convergence of $\theta_{N,h}$, which in turn allows us to produce consistent approximations of the optimal solution of the relaxed problem in eq.~\eqref{eq:original2}, as defined by Polak in~\cite[Sec.~3.3]{polak2012optimization}.
\begin{thm}
  $\theta_{N,h}$ epi-converges to $\theta_h$ as $N \to \infty$.
\end{thm}

The theorem above follows by Prop.~\ref{prop:sampling_converge}, Thm.~\ref{thm:approx}, and~\cite[Thm.~3.3.2]{polak2012optimization}.

\subsection{Input Signal Synthesis}
\label{subsec:synthesis}

Let $\mu_N^* \in {\cal M}_p(\R^m)$ be an empirical stochastic process built from a fixed set of samples $\set{u_i}_{i=1}^N \subset U$, as explained in Sec.~\ref{subsec:sampling}.
We synthesize a deterministic input signal following the procedure detailed in~\cite{vasudevan2013a,vasudevan2013b}.
Intuitively, given a sampling period $\Delta > 0$, the synthesis procedure generates an input $u^* \in L^2([0,T],\R^m)$ such that, for each integer $k < \frac{T}{\Delta}$:
\begin{equation}
  \int_{\mathrlap{k\Delta}}^{\mathrlap{(k+1)\Delta}} f\pb{t, x^{(u^*)}(t), u^*(t)}\, \diff{t} \approx
  \int_{\mathrlap{\R^m}} f\pb{t, x^{(\mu_N^*)}(t), u}\, \diff{\mu_{N,k\Delta}^*}(u).
\end{equation}

Denoting $\mu_{N,t}^* = \sum_{i=1}^N w_i^*(t)\, \delta_{u_i}$, the procedure aims to find a new set $\set{\hat{w}_i}_{i=1}^N$ whose values are strictly binary, i.e., $\hat{w}_i(t) \in \set{0,1}$ for each $t$.
We achieve this objective by first applying a Haar wavelet filter to each $w_i$, thus reducing them to piecewise constant functions.
Then, a pulse-width modulation (PWM) transform is applied to each filtered $w_i$, resulting in binary pulses whose widths are proportional to the amplitude of each function at the samples induced by~$\Delta$.
The procedure is explained in detail in~\cite[Sec.~4.4]{vasudevan2013a}.
Moreover, as shown in~\cite[Thm.~5.10]{vasudevan2013a}, $x^{(u^*)}$ converges to $x^{(\mu_N^*)}$ with rate~$\Delta^{1/2}$.

\begin{figure}[tp]
  \begin{algorithmic}[1]
    \Require $\set{u_i}_{i=1}^N \subset U$, $\set{w_i}_{i=1}^N \subset L^2([0,T],[0,1])$ with $\sum_{i=1}^N w_i(t) = 1$, $\epsilon_{\text{tol}} > 0$, $\alpha, \beta \in (0, 1)$.
    \State Let $\mu_{N,t} = \sum_{i=1}^N w_i(t)\, \one_{u_i}$.
    \Loop
    \State Compute $x^{(\mu_N)}$.
    \State Compute $\theta_{N,h}\pb{x^{(\mu_N)}, \mu_N}$ and $\delta \mu_N$ as in~\eqref{eq:ham_optimal_muN}.
    \If{$\theta_{N,h}\pb{x^{(\mu_N)}, \mu_N} > -\epsilon_{\text{tol}}$}
    \State Go to line~\ref{alg:numerical_break}.
    \EndIf
    \State Compute $\lambda^*$ as in~\eqref{eq:arimijo_mu}.
    \State $\mu_N \gets \mu_N + \lambda^*\, \delta \mu_N$.
    \EndLoop
    \State \label{alg:numerical_break}%
    Construct $\hat{w}_i \in L^2\pb{[0,T],\set{0,1}}$ equal to the PWM transform of a filtered $w_i$ using Haar wavelets, as explained in~\cite[Sec.~4.4]{vasudevan2013a}.
    \State \Return $u = \sum_{i=1}^N \hat{w}_i\, \one_{u_i}$.
  \end{algorithmic}
  \caption{Numerical algorithm to problem~\eqref{eq:original}.}
  \label{alg:numerical}
\end{figure}

Fig.~\ref{alg:numerical} summarizes our implementable algorithm.
Note that, in practice, we add an $\ell_1$ regularization term to the objective function of $\theta_{N,h}$, which results in sparse updates to $\mu_{N,t}$ on each iteration.
The sparse updates significantly improve the computation speed of the algorithm.
Moreover, since the $\ell_1$ regularization is strictly convex, it is not hard to show that this modification also results in an optimality function for the problem in eq.~\eqref{eq:original2}.

\section{SIMULATION RESULTS}
\label{sec:examples}


\subsection{Constrained LQR}
\label{subsec:example_1}

The first simulation is an input-constrained LQR problem. 
The system has 6~states, denoted $x$, and 2~inputs, denoted $u$.
The initial condition is $x(0) = 0$, and the vector field is:
\begin{equation}
  \smat{\ddot{x}_1\\ \ddot{x}_2\\ \ddot{x}_3} =
  \smat{-\sfrac{d}{m} & 0 & -\gamma\\ 0 & -\sfrac{d}{m} & 0\\ 0 & 0 & -\sfrac{m\, g\, l}{J}}\,
  \smat{\dot{x}_1\\ \dot{x}_2\\ x_3} +
  \smat{\sfrac{1}{m} & 0\\ 0 & \sfrac{1}{m}\\ \sfrac{r}{J} & 0}\,
  \smat{u_1\\ u_2},
\end{equation}
with objective function:
\begin{equation}
  \int_0^T \norm{\smats{x_1(t) + 0.3\\ x_2(t) + 0.5\\ x_3(t)}}_2^2 + \eta\, \norm{u(t)}_2^2\, \diff{t}.
\end{equation}
The parameters are $T=2$, $J = 0.0475$, $m = 1.5$, $r = 0.25$, $g = 9.8$, $\gamma = 0.51$, $d = 0.2$, $l = 0.05$, and $\eta = 0.05$.
We use 81~points to sample the vector field, evenly distributed in the control space $U = [-1,1]^2$.

\begin{figure*}
  \centering
  \subfloat[%
    Optimal trajectory.
    Final state shown as a triangle, $x(T) = (-0.29, -0.56, -0.06)$.%
  ]{%
    \label{fig:exp1_x}%
    \includegraphics[width=.31\linewidth,trim=20 0 0 10,clip]{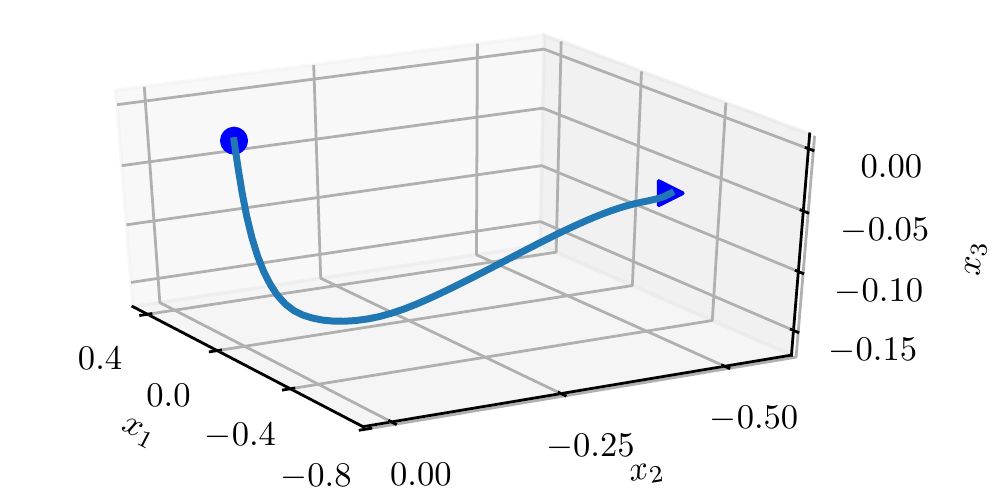}%
  }%
  \hfill%
  \subfloat[%
    PWM projection of the optimal inputs.
    Both inputs are identically zero for each $t \geq 0.5$.
  ]{%
    \label{fig:exp1_u}%
    \includegraphics[width=.31\linewidth,trim=0 0 0 0,clip]{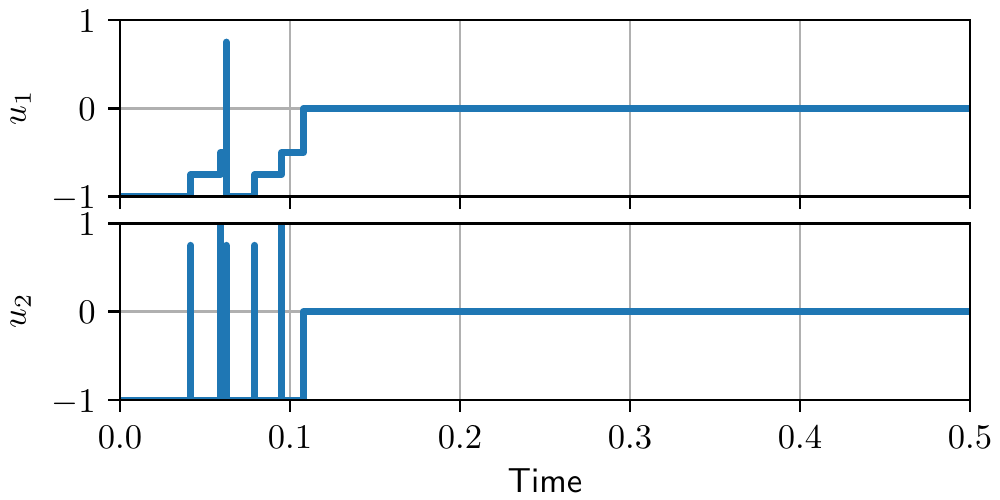}%
  }%
  \hfill%
  \subfloat[%
    Optimality function value per iteration.
  ]{%
    \label{fig:exp1_theta}%
    \includegraphics[width=.31\linewidth,trim=0 0 0 0,clip]{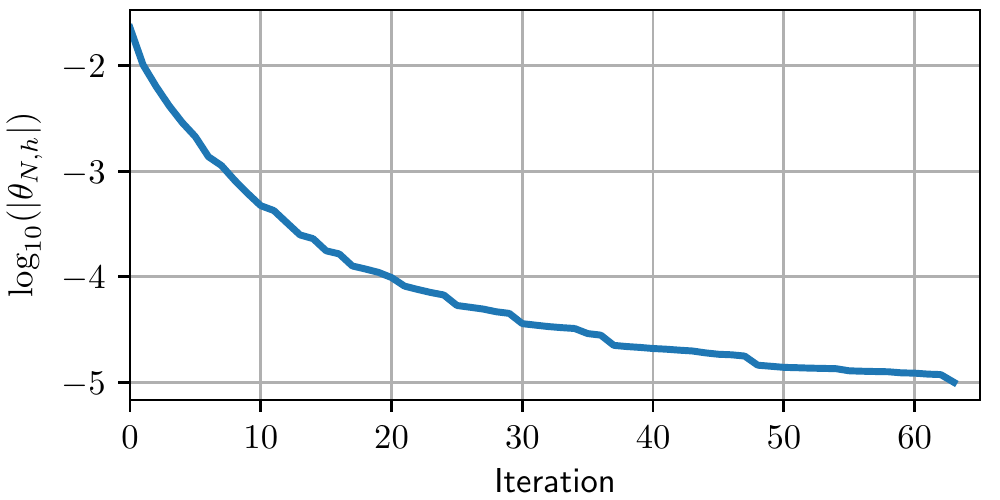}%
  }%
  \caption{%
    Results of the constrained LQR simulation in Sec.~\ref{subsec:example_1}.
  }
  \label{fig:exp1}
\end{figure*}

Fig.~\ref{fig:exp1_x} shows the trajectory of the first three states of the system. 
Fig.~\ref{fig:exp1_u} shows the corresponding controls after wavelet and PWM reconstruction.
Note how the input signal effort is very low, only applying an input for the first $0.1\unit{s}$ approx.
Fig.~\ref{fig:exp1_theta} shows the absolute value of the optimality function $\theta_{N,h}$ at each iteration.
After 62~iterations the optimality function reaches the threshold $\epsilon_{\text{tol}} = 10^{-5}$, validating the convergence of our algorithm.

\subsection{Quadrotor helicopter}
\label{subsec:example_2}

We also consider an application to control a 3-D quadrotor helicopter using a nonlinear model described in~\cite{pounds2006modelling}.
This system has 12~states, denoted $x$, and 4~inputs, denoted $u$.
The inputs are constrained to the set $U = [0, 2]^4$.
The initial condition is $x(0) = 0$, and the vector field is defined as follows:
\begin{equation}
  \begin{aligned}
    \ddot{x}_1 &= - \sfrac{b}{m}\, \dot{x}_1 + \sfrac{K}{m}\, \sin(x_5)\, \smats{1 & 1 & 1 & 1}\, u,\\
    \ddot{x}_2 &= - \sfrac{b}{m}\, \dot{x}_2 + \sfrac{K}{m}\, \sin(x_4)\, \cos(x_5)\, \smats{1 & 1 & 1 & 1}\, u,\\
    \ddot{x}_3 &= - \sfrac{b}{m}\, \dot{x}_3 + \sfrac{K}{m}\, \cos(x_4)\, \cos(x_5)\, \smats{1 & 1 & 1 & 1}\, u,\\
    \ddot{x}_4 &= - \dot{x}_5\, \dot{x}_6 + \sfrac{K\, L}{I_x}\, \smats{0 & 1 & 0 & -1}\, u,\\
    \ddot{x}_5 &= \dot{x}_5\, \dot{x}_6 + \sfrac{K\, L}{I_y}\, \smats{-1 & 0 & 1 & 0}\, u,\\
    \ddot{x}_6 &= \sfrac{K}{I_x + I_y}\, \smats{1 & -1 & 1 & -1}\, u,
  \end{aligned}
\end{equation}
with objective function:
\begin{equation}
  \int_0^T \norm{\smats{x_1(t) - c_1\\ x_2(t) - c_2\\ x_3(t) - c_3}}_2^2 + \sum_{k=4}^6 \sin^2\pb{x_k(t)} + \eta\, \norm{u(t)}\, \diff{t}.
\end{equation}
The states $x_{1, 2, 3}$ represent the $x$, $y$, and $z$ position coordinates, respectively.
The states $x_{4, 5, 6}$ represent the Euler angles of the body.
The parameters are $m = 1.3 \unit{kg}$, $I_{x, y}=0.0605\unit{kg\cdot m^2}$, $g = 9.8\unit{m/s^2}$, $b=0.1$, $K = 1.0$, and $\eta = 0.05$.
We use 625~sampling points for the vector field, and the convergence threshold is $\epsilon_{\text{tol}} = 10^{-4}$.

\begin{figure}[tp]
  \centering
  \subfloat[%
    $p_f = (-1.12, -1.12, -0.98)$.%
  ]{%
    \label{fig:exp2_1}%
    \includegraphics[width=.48\linewidth,trim=15 0 0 10,clip]{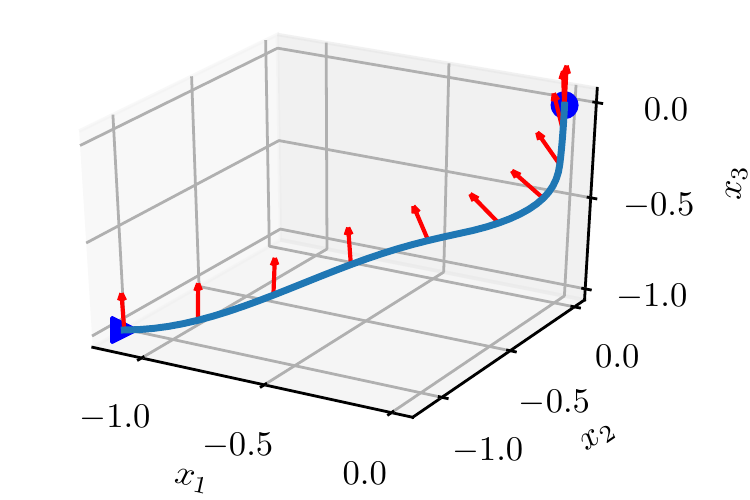}%
  }%
  \hfill%
  \subfloat[%
    $p_f = (1.23, 1.23, 1.18)$.%
  ]{%
    \label{fig:exp2_2}%
    \includegraphics[width=.48\linewidth,trim=15 0 0 10,clip]{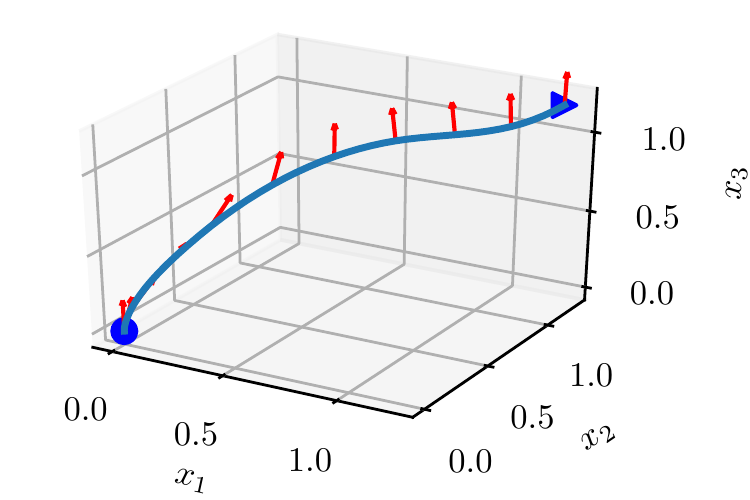}%
  }%
  \caption{%
    Results of the quadrotor helicopter simulation in Sec.~\ref{subsec:example_2}.
    Final position, $p_f = \sp{x_1(T), x_2(T), x_3(T)}$, is shown as a triangle, and vertical orientation of the helicopter is shown as red arrows.
  }
  \label{fig:exp2}
\end{figure}

\begin{figure}[tp]
  \centering
  \includegraphics[width=.67\linewidth,trim=0 0 0 0,clip]{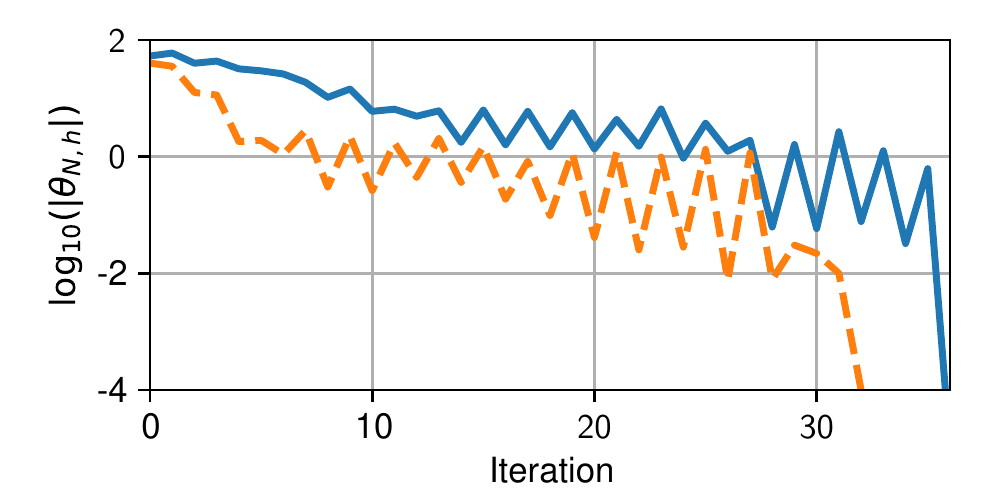}%
  \caption{%
    Optimality function value per iteration for the simulations in Fig.~\ref{fig:exp2_1} (blue solid line) and Fig.~\ref{fig:exp2_2} (orange dashed line).
  }
  \label{fig:exp2_3}%
\end{figure}

Fig.~\ref{fig:exp2} shows the results of the optimal trajectories for two scenarios: when $c = (-1.2, -1, -1)$ (Fig.~\ref{fig:exp2_1}), and when $c = (1.2, 1, 1)$ (Fig.~\ref{fig:exp2_2}).
Fig.~\ref{fig:exp2_3} shows the absolute value of the optimality function $\theta_{N,h}$ at each iteration in those two simulations.
In both cases we reach the convergence threshold $\epsilon_{\text{tol}}$ in less than 40~iterations.
These results show that our algorithm generates trajectories that move the helicopter to the desired positions while at the same time resulting in minor body-axis displacements.

\section{CONCLUSION}
\label{sec:conclusion}

In this paper we present a theoretical formulation, and a corresponding numerical algorithm that can find Pontryagin-optimal inputs for general dynamical systems by using a direct method.
The numerical implementation is based on a relaxed-control system and PWM reconstruction.
Our novel approach produces significant improvements both in the quality of the resulting minimizers, and the flexibility of the numerical implementation.




\bibliographystyle{IEEEtran}
\bibliography{refs}

\end{document}